\documentclass[11pt]{article}
\usepackage {amsfonts}
\usepackage[english]{babel}
\begin{document}

\large

\def\Subset{\subset\subset}
\def\card{{\rm{card}}}
\def\C{{\mathbb{C}}}
\def\R{{\mathbb{R}}}
\def\Z{{\mathbb{Z}}}
\def\N{{\mathbb{N}}}
\def\Im{{\rm{Im}}}
\def\Re{{\rm{Re}}}
\def\s{{\sigma}}
\def\d{\delta}
\def\r{\rho}
\def\f{\varphi}
\def\a{\alpha}
\def\b{\beta}
\def\t{\tau}
\def\th{\theta}
\def\e{{\varepsilon}}
\def\l{\lambda}
\def\g{\gamma}
\def\Arg{{\arg}}
\def\supp{{\rm{supp}}}

\newtheorem {theor}{Theorem}
\newtheorem {lem}{Lemma}
\newcommand {\theorFL}{Theorem}
\newcommand {\bproof} {\noindent{\bf Proof} }
\newcommand {\eproof} {\noindent}

\title{
    \bf Subharmonic Almost Periodic Functions of Slowly Growth}

\author{Favorov S.Ju., Rakhnin A.V.}

\date{}

\maketitle

\begin{abstract}
We obtain a complete description of the Riesz measures of almost
periodic subharmonic functions with at most of linear growth on
$\C$; as a consequence we get a complete description of zero sets
for the class of entire functions of exponential type with almost
periodic modulus.
\end{abstract}



\bigskip

Bohr's Theorem (see \cite{B}, or \cite{l}, Ch.6, \S 1) implies
that each almost periodic function\footnote{Explicit definitions
of almost periodicity for functions, measures, and discrete sets
see \cite{l}, Ch.6 and Appendix VI, \cite{r}, or \S 3 of the
present paper.} on real axis $\R$ with the bounded spectrum is
just the restriction to $\R$ of an entire almost periodic function
$f$ of exponential type. Moreover, $f$ has no zeros outside of
some strip $|\Im z| \le H$ if and only if supremum and infimum of
the spectrum $f$ also belongs to the spectrum. In \cite{kl} (see
also \cite{l}, Appendix VI) M.G.Krein and B.Ya.Levin obtained a
complete description of zeros of functions from the last class.
Namely, a set $\{a_k\}_{k\in\Z}$ in a horizontal strip of a finite
width is just the zero set of an entire almost periodic function
$f$ of exponential type if and only if the set is almost periodic
and has the representation
\begin{equation}\label{ak}
    a_k = d k + \psi (k), \; k\in\Z,
\end{equation}
where $d$ is a constant, the function $\psi (k)$ is bounded, and
the values
\begin{equation}\label{sn}
S_n = \lim\limits_{r\to\infty} \sum\limits_{|k|<r} [\psi (k+n) -
\psi (k)]\frac{k}{k^2+1}
\end{equation}
are bounded uniformly in $n\in\Z$.

It can be proved that almost periodicity of $\{a_k\}$ yields
representation (\ref{ak}) and a finite limit in (\ref{sn}) for
every fixed $n\in\Z$. Also, one can obtained a complete
description of zero sets for the class of entire functions of
exponential type with almost periodic modulus and zeros in
horizontal strip of finite width: we should only replace the $S_n$
by $\Re S_n.$

Observe that every entire function of exponential type bounded on
$\R$ has the form
\begin{equation}\label{predst-f}
    f(z) = C e^{i\nu z} \lim\limits_{r\to\infty}
    \prod\limits_{|a_k|<r}\left(1-\frac{z}{a_k}\right), \;\nu \in
    \R
\end{equation}
(\cite{l}, Ch. 5; for simplicity we suppose $0\not\in \{a_k\}$);
so we have an explicit representation for functions from the
classes above-mentioned.

Note that one of the authors of the present paper obtained in
\cite{f1} a complete description of zero sets for holomorphic
almost periodic functions on a strip and on the plane without any
growth conditions. An implicit representation for a special case
of almost periodic holomorphic functions was obtained earlier in
\cite{rrf}. Besides, it was proved in \cite{rrf} that zero sets of
holomorphic  functions with the almost periodic modulus on a strip
(or on the plane) are just almost periodic discrete sets. This
result is a consequence of more general one: every almost periodic
measure on a strip is just the Riesz measure of some subharmonic
almost periodic function on the strip.

In \S 3 of our paper we obtain a complete description of the Riesz
measures for almost periodic subharmonic functions of the normal
type with respect to the order 1 (note that it is the smallest
growth for bounded on $\R$ subharmonic function). In particular,
we consider the case of periodic subharmonic functions. As a
consequence, we get a complete description of zero sets for the
class of entire functions of exponential type with the almost
periodic modulus without any additional requirements on
distributions of zeros. Note that representation (\ref{ak}) with a
bounded function $\psi (k)$ is incorrect here, therefore methods
of paper \cite{kl} do not work in our case. The integral
representation from \cite{rrf} creates an almost periodic
subharmonic function with a given almost periodic Riesz measure,
but does not allow to control the growth of the function,
therefore it isn't fit for our problem as well.

We make use of a subharmonic analogue of representation
\begin{equation}\label{predst-log}
    \log |f(z)| = \int\limits_{0}^{\infty} \frac{n(0,t) -
    n(z,t)}{t}dt - \nu y + \log |C|,
\end{equation}
for functions of the form (\ref{predst-f}) (see \cite{f2} or
review \cite{Kh}, p.45); here $n(c,t)$ is a number of zeros in the
disc $\{ z: \; |z-c|\le t \}$. We obtain this analogue in \S 2 of
our paper. Also, we get a complete description of the Riesz
measures for bounded on $\R$ subharmonic functions with at most of
linear growth on $\C$.

Here we base on a subharmonic analogue of (\ref{predst-f}) as
well. Of course, this analogue can be obtained by repeating all
the steps of the proof (\ref{predst-f}) for entire functions in
\cite{l}, nevertheless we prefer to give a short proof in \S 1,
using Azarin's theory of limit sets for subharmonic functions
\cite{a}. The idea of the proof belongs to prof. A.F.Grishin, and
the authors is very grateful to him.

\bigskip
\begin{center}
{\bf \S 1}
\end{center}
\bigskip

In this section we prove the following theorem:

\begin{theor}.
Let $v(z)$ be a subharmonic function on $\C$ such that

\begin{equation}\label{4}
    v^{+}(z) = O(|z|)\;\hbox{as} \; z\to\infty
\end{equation}
and
\begin{equation}\label{5}
    \sup\limits_{x\in\R}v(x) < \infty.
\end{equation}
Then
\begin{equation}\label{5a}
    v(z) = \lim\limits_{R\to\infty} \int\limits_{|w|<R} (\log|z-w| - \log^{+}
    |w|)d\mu(w) + A_1 y + A_2.
\end{equation}
Here $z=x+iy$, $\mu = \frac{1}{2\pi} \triangle v$ is the Riesz
measure of $v$, $A_1, A_2 \in \R$, and the limit exists uniformly
on compact subsets in $\C$.
\end{theor}

Note that the condition (\ref{4}) means just $v$ is at most of
normal type with respect to the order 1.

Our proof of Theorem 1 based on Azarin's theory of limit sets
\cite{a}. Thus, if a subharmonic function $v$ satisfies (\ref{4}),
then

a) the family $v_t(z) = t^{-1} v(tz)$, $t>1$ is a relatively
compact set in the space of distributions $\cal D'(\C)$; in other
words, for every sequence of functions from this family there is a
subsequence converging to a subharmonic function.

Note that the convergence in $\cal D'(\C)$ is the weak convergence
on all functions from the class of infinitely smooth compactly
supported functions $\cal D(\C)$; moreover, the class of
subharmonic functions is closed with respect to this convergence.

b) If
\begin{equation}\label{6}
    v_{\infty} = \lim\limits_{t\to\infty} v_t(z),
\end{equation}
then the Riesz measure $\mu_{\infty}$ of the function $v_{\infty}$
satisfies the equality
\begin{equation}\label{7}
    \mu_{\infty} = \lim\limits_{t\to\infty} \mu_t,
\end{equation}
where $\mu_t(E) = t^{-1} \mu (tE)$ for Borel subsets of $\C$;
limits (\ref{7}) exists in the sense of weak convergence on
continuous compactly supported functions on $\C$. Moreover, in
this case there exists
\begin{equation}\label{8}
    \lim\limits_{R\to\infty} \int\limits_{1\le|w|<R}
    \frac{d\mu(w)}{w}\ne \infty.
\end{equation}

If a subharmonic function satisfies (\ref{4}) and (\ref{6}), then
it is called {\it a completely regular growth (with respect to the
order 1)}.

\bigskip
In what follows we need the simple criterion of compactness for
family of subharmonic functions (see\cite{a}).

\medskip

{\bf Lemma A}. A family $\{u_\a\}$ of subharmonic functions on
$\C$ is a relatively compact set in the space of distributions
$\cal D'(\C)$ if and only if

a) $\sup_\a\sup_{z\in K}u_\a(z)<\infty$ for all compacta
$K\subset\C$,

b) $\inf_\a\sup_{z\in K_0}u_\a(z)>-\infty$ for some compact set
$K_0\subset\C$.

\bigskip

Also, we need the following variant of Fragmen-Lindelof Theorem.

\medskip

{\bf Theorem FL}. {\it If a function $v$ is subharmonic in a
neighborhood of the closure of upper half-plain  $\C^{+} =
\{z=x+iy:\; y>0\}$ and satisfies conditions (\ref{4}),  (\ref{5}),
then  for all $z\in\C^{+}$
$$
v(z) \le \sup\limits_{x\in\R} v(x) + \s^+ y,
$$
with $\s^+ = \limsup\limits_{y\to +\infty}y^{-1} v(iy)$.}

The proof of this statement is the same as for holomorphic on
$\C^{+}$ and continuous on $\overline{\C^{+}}$ functions (see, for
example, \cite{k}, p.28).

\medskip

\bigskip

First let us prove a subharmonic analogue of Cartwright Theorem
(the holomorphic case see, for example, \cite{l}, Ch. V).

\begin{theor}.
Let a subharmonic function $v$ on $\C$ be satisfied (\ref{4}) and
(\ref{5}). By definition, put
\begin{equation}\label{11}
    \s_{\pm} = \limsup\limits_{y\to\pm\infty} \frac{v(iy)}{|y|}.
\end{equation}
Then $v(z)$ is a completely regular growth; moreover, the function
$v_{\infty}$ from (\ref{6}) has the form
\begin{eqnarray} \label{10a}
v_{\infty}(z) = \left\{
\begin{array}{rl}
    \s^+ y, & y\ge0, \\
    \s^- |y|, & y<0, \\
\end{array}
\right.
\end{eqnarray}
\end{theor}

\medskip

{\bf Remark}. From Theorem FL it follows that if a subharmonic
function $v$ on $\C$ satisfies the conditions of Theorem 2 with
$\s^+\le0$ and $\s^-\le0$, then $v$ is a constant.

\medskip

The proof is based on the following lemma

\begin{lem}.
Let $u<0$ be a subharmonic function on $\C^{+}$. Then for every
$R<\infty$ and $r\in(0,R/2)$
\begin{equation}\label{9}
    \frac{u(ir)}{r} \le \frac{C}{R^3}
    \int\limits_{|z-iR|<R} u(z)dm_2(z),
\end{equation}
where $C$ is an absolute constant and  $m_2$ is the plain Lebesque
measure.
\end{lem}

\bproof

 Since  Poisson formula for the disc $B(iR/2,R/2) =
\{z:\; |z-iR/2|<R/2\}$, we have

$$
u(ir) \le \frac{1}{2\pi} \int \limits_{0}^{2\pi} u(iR/2 +
e^{i\theta} R/2) \frac{(R/2)^2 - (R/2 - r)^2}{(R/2)^2 + (R/2 -
r)^2 - R(R/2 - r)\cos (\pi/2 + \theta)}d\theta .
$$
 Using the inequality $u<0$, replace the interval of integration by
$[\pi/4, 3\pi/4]$. We obtain
\begin{equation}\label{10}
    u(ir) \le \frac{r}{(R-r)} \frac{1}{4} \sup\limits_{\theta \in [\pi/4, 3\pi/4]}
    u \left(\frac{iR}{2} + \frac{R}{2}e^{i\theta}\right)
\end{equation}
Since for all $\theta \in [\pi/4, 3\pi/4]$
$$
B(iR,R/2) \subset B(iR/2 + e^{i\theta} R/2, (1 +
\frac{\sqrt{2}}{2} )R/2) \subset \C^{+},
$$
we have
$$
u(\frac{iR}{2} + \frac{R}{2} e^{i\theta}) \le \frac{4}{\pi R^2
(1+\frac{\sqrt{2}}{2})^2} \int\limits_{|z-iR/2-e^{i\theta}R/2|\le
(1+\frac{\sqrt{2}}{2})R/2} u(z)dm_2(z)\le
$$
$$
\le\frac{8}{\pi R^2 (3+2\sqrt{2})} \int\limits_{|z-iR|<R/2}
u(z)dm_2(z)
$$
Then replace the average over the disc $B(iR,R/2)$ by the average
over the disc $B(iR,R)$. So the assertion of the Lemma follows
from (\ref{10}).

 \eproof

\medskip

{\bf Proof of theorem~2}.

Without loss of generality it can be assumed that
$\sup\limits_{x\in \R} v(x) = 0.$ Put $u(z) = v(z) - \s_{+} y$.
From Theorem FL it follows that $u(z) < 0$ on $\C^{+}$, then
\begin{equation}\label{12}
     \limsup\limits_{y\to +\infty} \frac{u(iy)}{y} = 0
\end{equation}
Fix $z_0\in \C^{+}$. Let $\f$, $0\le\f\le1$, be an infinite
differentiable and compactly supported function on $\C^{+}$,
depending only on $|z-z_0|$. Apply Lemma 1 for the function
$u_t(z) = u(tz) t^{-1}$. If $\supp \f \subset B(iR,R)$, then we
get
$$
\frac{u(itr)}{tr}\le \frac{C}{R^3}
\int\limits_{|z-iR|<R} u_t(z)dm_2(z) \le \frac{C}{R^3}
\int\limits_{\C^{+}} u_t(z) \f(z) dm_2(z)
$$
Since (\ref{12}), we see that for each $\e>0$ and $t>t(\e)$ there
is $r\in (0,R/2)$ such that $u(itr) \ge -\e tr$. We obtain
$$
\int\limits_{\C^{+}} u_t(z) \f(z) dm_2(z) \ge -\frac{\e R^3}{C}.
$$
Hence,
$$
\int\limits_{\C^{+}} u_t(z) \f(z) dm_2(z) \to 0 \quad \hbox{as}
\quad t\to\infty.
$$
Therefore, $v_t(z)=u_t(z)+\s_+ y \to \s_{+} y$ in the space $\cal
D'(\C^{+})$. Similarly, $v_t(z) \to -\s_{-} y$ in the space $\cal
D'(\C^{-})$, where $\C^{-} = \{z=x+iy: \; y<0\}$. Each limit
function for $v_t(z)$ is always subharmonic, therefore we get
(\ref{10a}). So limit (\ref{6}) exists. Theorem is proved.
 \eproof

\bigskip

{\bf Consequence.} {\it The Riesz measure of the limit function
$v_{\infty}(z)$ equals
$$
\frac{\s_{+} + \s_{-}}{2\pi} m_1(x),
$$
where $m_1$ being the Lebesgue measure on $\R$. }

{\bf Proof of theorem~1}. From Jensen-Privalov formula for
subharmonic function, we get the estimate for $r\ge 1$
\begin{equation}\label{13}
    \mu (B(0,r))\le C_1 r,
\end{equation}
where constant $C_1$ depends only on $v$. Using Brelot-Hadamard
Theorem for subharmonic function (see, for example, \cite{r1}), we
obtain that there is  a harmonic polynomial $H(z)$ of degree 1
such that
\begin{equation}\label{14}
    v(z) = \int\limits_{|w|<1}\log |z-w|d\mu(w) + \int\limits_{|w|\ge
    1}\log \left(\left|1-\frac{z}{w}\right| + \Re
    \frac{z}{w}\right)d\mu(w) + H(z).
\end{equation}
Denote by $v^0(z)$ the first integral in (\ref{14}). Since
(\ref{8}), we get
\begin{equation}\label{15}
    v(z) = v^0(z) + \lim\limits_{R\to\infty} \int\limits_{1\le
    |w|<R}\log \left|1-\frac{z}{w}\right|d\mu(w) + A_0 x + A_1 y +
    A_2.
\end{equation}
The application of Theorem 2 yields that the function $v$ is a
completely regular growth, hence the measures $\mu_t$ converge
weakly to the measure
\begin{equation}\label{15a}
\mu_{\infty} = \frac{\s_{+} + \s_{-}}{2\pi} m_1.
\end{equation}
Let $\mu'$ be the restriction of the measure $\mu$ to $\C\setminus
B(0,1)$. Obviously, measures $\mu'_t$ weakly converge to the
measure $\mu_{\infty}$ as well. Therefore,
\begin{equation}\label{16}
    v_t(z) = v^0_t(z) + \lim\limits_{R\to\infty}
    \int\limits_{|w|<R}\log \left|1-\frac{z}{w}\right|d\mu_{t}'(w) + A_0 x + A_1 y +
    \frac{A_2}{t}.
\end{equation}

Pass to a limit in (\ref{16}) as $t\to\infty$ in the space $\cal
D'(\C)$. First, by Theorem 2, the functions $v_t(z)$ converge to
the function $v_{\infty}(z)$ from (\ref{10a}). Since $v^0(z) =
O(\log|z|)$ as $|z|\to\infty$, we see that $v^0_t(z)\to 0$. By
(\ref{13}), we obtain $\mu'_t (B(0,R))\le C_1 R$ for all $t\ge 1$
and $R>0$. Therefore, we get uniformly in $t \ge 1$,
\begin{equation}\label{17}
    \int\limits_{|w| \ge R} \frac{d\mu_{t}'(w)}{|w|^2} = 2\int\limits_{R}^{\infty}
    \frac{\mu_{t}(B(0,s))}{s^3}ds - \frac{{\mu_t}'(R)}{R^2} \to 0.
\end{equation}
as $R\to\infty$. Also, by (\ref{8}), uniformly in $t\ge 1$, $R'\ge
R$
\begin{equation}\label{18}
\int\limits_{R \le |w| \le R'} \frac{d\mu_{t}'(w)}{w} =
\int\limits_{Rt \le |w| \le R't} \frac{d\mu(w)}{w} \to 0,
\end{equation}
as $R\to\infty$. For $|z|<C$ and a sufficiently large $|w|$
$$
\left|\log\left|1-\frac{z}{w}\right| + \Re \frac{z}{w}\right|\le
\frac{|z|^2}{|w|^2}.
$$
Therefore, taking into account (\ref{17}) and (\ref{18}), we
obtain for all $\f\in\cal D(\C)$ uniformly in $t \ge 1$ $$
    \int\limits_{\C} \left( \lim\limits_{R\to\infty} \int\limits_{|w|\le R} \log
    \left|1-\frac{z}{w}\right|d{\mu_t}'(w)\right)\f (z) dm_2(z) =
$$
\begin{equation}\label{19}
    = \lim\limits_{R\to\infty} \int\limits_{|w|\le R}\left( \int\limits_{\C}\log
    \left|1-\frac{z}{w}\right|\f (z)dm_2(z) \right)d{\mu_t}'(w).
\end{equation}
Note that measure $\mu_{\infty}$ does not charge any circle
$|w|=R$, therefore the restrictions of measures ${\mu_t}'(w)$ to
any disc $B(0,R)$ weakly converge to the restriction of the
measure $\mu_{\infty}(w)$. The function $ \int\log|z-w|\f (z)
dm_2(z)$ is continuous in the variable $w$, so we have
\begin{equation}\label{20}
    \lim\limits_{t\to\infty} \int\limits_{|w| \le R} \int\log|z-w|\f (z)
    dm_2(z)d\mu_t(w) =  \int\limits_{|w| \le R} \int\log|z-w|\f (z)
    dm_2(z)d\mu_{\infty}(w).
\end{equation}
By the same reason for each $\d> 0$
\begin{equation}\label{21}
    \lim\limits_{t\to\infty} \int\limits_{\d \le |w| \le R}
    \log|w|\mu_t(w) = \int\limits_{\d \le |w| \le R}
    \log|w|\mu_{\infty}(w).
\end{equation}
Furthermore,
\begin{equation}\label{24a}
\int\limits_{|w| \le \d} \log|w|{\mu'}_t(w) = \log \d \;
\frac{\mu'(B(0,\d t))}{t} - \int\limits_{0}^{\d}
\frac{\mu'(B(0,st))}{st}ds
\end{equation}

Since $\mu'(B(0,r))\le C_1 r$ for all $r>0$, we see that
(\ref{24a}) tends to zero as $\d\to 0$ uniformly in $t\ge 1$.
Combining (\ref{19}), (\ref{20}), (\ref{21}), and (\ref{24a}), we
get the equality
$$
v_{\infty} (z) = \lim\limits_{R\to\infty}\int\limits_{|w|\le R}
\log\left|1-\frac{z}{w}\right|d\mu_{\infty}(w) + A_0 x + A_1 y.
$$
Take $y=0$. Since $v_{\infty} (x) = 0$ and
$$
\lim\limits_{R\to\infty}\int\limits_{|u|\le R}
\log\left|1-\frac{x}{u}\right|du = 0
$$
for all $x\in\R$, we obtain $A_0 = 0$. Now the assertion of
Theorem 1 follows from (\ref{15}).

\bigskip

\bigskip
\begin{center}
{\bf \S 2}
\end{center}
\bigskip

Here we get a complete description of Riesz measures for
subharmonic functions with at most of linear growth on $\C$ (i.e.,
do not exceed $C(|z|+1)$ with $C<\infty$) and with some additional
conditions (bounded on $\R$ or with the compact family of
translations along $\R$). Holomorphic analogues of the
corresponding theorems were obtained earlier one of the author in
\cite{f2}.

First prove some lemmas.

\begin{lem}. If a measure $\mu$ on $\C$ satisfies the condition
\begin{equation}\label{22}
    \mu(B(0,R+1)) - \mu(B(0,R)) = \overline{o}(R) \; \hbox{as} \;
    R\to\infty,
\end{equation}
and the limit
$$
\lim\limits_{R\to\infty} \int\limits_{|w|<R} (\log^+ |z-w| -
\log^+ |w|)d\mu(w)
$$
exists at some point  $z\in\C$, then the limit equals
$$
\int\limits_{1}^{\infty} \frac{\mu(B(0,t)) - \mu(B(z,t))}{t}dt.
$$
\end{lem}

\bproof

For all $z\in\C$ and $R\in(|z|+1,\infty)$ we have
$$
\int\limits_{|w|<R}\log^+ |z-w| d\mu(w) -
\int\limits_{|w|<R}\log^+ |w| d\mu(w) = (\log R) \mu (B(z,R)) -
$$
$$
-\int\limits_{1}^{R}\frac{\mu(B(z,t))}{t}dt - (\log R) \mu(B(0,R))
+ \int\limits_{1}^{R}\frac{\mu(B(0,t))}{t}dt +
$$
$$
 + \int\limits_{|w|<R, \;
|w-z|\ge R}\log^+ |z-w| d\mu(w) - \int\limits_{|w|\ge R, \; |w-z|<
R}\log^+ |z-w| d\mu(w) =
$$
$$
=\int\limits_{1}^{R}\frac{\mu(B(0,t)) - \mu(B(z,t))}{t}dt +
$$
\begin{equation}\label{22a}
+ \int\limits_{|w|<R, \; |w-z|\ge R}\log
\left|\frac{z-w}{R}\right| d\mu(w) - \int\limits_{|w|\ge R, \;
|w-z|< R}\log \left|\frac{z-w}{R}\right| d\mu(w).
\end{equation}
If $|w|<R$ and $|z-w| \ge R$ or $|w|\ge R$ and $|z-w| < R$, then
we have
$$
1 - \frac{|z|}{R}\le \left| \frac{z-w}{R}\right| \le 1+
\frac{|z|}{R}.
$$
Therefore the integrand
functions of last two integrals in (\ref{22a}) are $O(1/R)$ as
$R\to\infty$. The domains of integrations are subsets of the ring
$R-|z| \le |w| \le R+ |z|$, hence, by (\ref{22}), these integrals
tends to $0$ as $R\to\infty$. Lemma is proved.

 \eproof

\begin{lem}.
Let a measure $\mu$ be satisfied (\ref{8}), (\ref{13}), and
(\ref{22}), and let
\begin{equation}\label{23}
    V(z) = \lim\limits_{R\to\infty} \int\limits_{|w|<R} (\log|z-w| - \log^{+}
    |w|)d\mu(w).
\end{equation}
Then $V$ is a subharmonic function with Riesz measure $\mu$ and
\begin{equation}\label{24b}
    V(z) = \int\limits_{1}^{\infty} \frac{\mu(B(0,t)) -
    \mu(B(z,t))}{t}dt+ \int\limits_{|w-z| < 1} \log |z-w|d\mu(w)
\end{equation}
for all $z\in\C$. Furthermore, the function
\begin{equation}\label{28}
    \tilde V(z) = \frac{1}{2\pi} \int\limits_{0}^{2\pi} V(z+e^{i\theta})
    d\theta
\end{equation}
satisfies the equality
\begin{equation}\label{24}
    \tilde V(z) =\int\limits_1^{\infty} \frac{\mu(B(0,t)) -
    \mu(B(z,t))}{t}dt.
\end{equation}
\end{lem}

\bproof

It follows from (\ref{13}) that the integral in (\ref{14}) is a
subharmonic function on $\C$ with the Riesz measure $\mu$;
besides, it satisfies (\ref{4}) (see, for example, \cite{r1},
Ch.1). If, in addition, $\mu$ satisfies (\ref{8}), then the limit
in (\ref{23}) exists uniformly on bounded sets, and the function
$V$ coincides with the integral in (\ref{14}) up to a linear term;
so it has the same properties as well.

 Using the equality $\int_{0}^{2\pi}
\log|a+e^{i\theta}|d\theta=2\pi\log^+|a|$, we get
$$
\tilde V(z) = \lim\limits_{R\to\infty} \int\limits_{|w|<R}
(\log^+|z-w| - \log^+ |w|)d\mu(w).
$$
and $V(z)=\tilde V(z)+\int_{|w-z| < 1} \log |z-w|d\mu(w)$. Then
Lemma 2 implies (\ref{24}) and (\ref{24b}). Lemma 3 is proved.
\eproof

\begin{lem}.
Let a subharmonic on $\C$ function $v$ be satisfied (\ref{4}).
Then the family of translations $\{v(z+h)\}_{h\in\R}$ is a
relatively compact subset in $\cal D'(\C)$ if and only if the
function
\begin{equation}\label{31a}
    \tilde v(z) = \frac{1}{2\pi} \int_{0}^{2\pi} v(z+e^{i\theta})d\theta
\end{equation}
is bounded on $\R$; this function is bounded simultaneously with
the function
\begin{equation}\label{31b}
    \hat v(z) = \frac{1}{\pi} \int_{|w|<1} v(z+w)dm_2(w)
\end{equation}

\end{lem}

\bproof. If the family $\{v(z+h)\}_{h\in\R}$ is a relatively
compact subset, then it is uniformly bounded from above on
compacta in $\C$ and the function $v$ is uniformly bounded from
above on every strip $|y|<H$. Then by Lemma A there is a compact
subset $K_0$ of $\C$ such that
$$
\sup \limits_{K_0} v(z+h) \ge C_2,\quad\forall h\in\R.
$$
Take $d>\sup_{K_0}|z|$. Then for each $h\in\R$ there is a point
$z(h)$, $|z(h)|<d$, such that
 $$
\frac{1}{(d+1)^2\pi} \int \limits_{|w-z(h)|\le d+1} v(w+h)dm_2(w)
\ge v(z(h) + h) > C_2 - 1.
 $$
 Further, we have
 $$
\int \limits_{|w-z(h)|\le d+1} v(w+h)dm_2(w)
\le\int\limits_{|w|\le 1} v(w+h)dm_2(w)
+(d^2+2d)\pi\sup_{|y|<d+1}v(z).
 $$
Therefore,  $\inf_{h\in\R} \hat v(h)> -\infty$. Since $\tilde
v(h)\ge\hat v(h)$, we see that the functions $\hat v$ and $\tilde
v$ are bounded uniformly from below on $\R$. It is clear that
these functions are bounded uniformly from above on $\R$ as well.

On the other hand, if $\tilde v(z)$ is bounded from below on $\R$,
then $\inf_{h\in\R}\sup_{|w|=1}v(h+w)>-\infty$; if $\tilde v(h)$
is bounded from above on $\R$, then $v(h)$ is bounded from above
on $\R$ as well. Since Theorem FL, we see that $v(z)$ is bounded
from above on every strip $|y|<H$. It follows from Lemma A that
$\{v(z+h)\}_{h\in\R}$ is a relatively compact set. Hence $\hat
v(x)$ is bounded on $\R$.\eproof
\bigskip

Now we can prove the theorems mentioned above.

\begin{theor}.
For a measure $\mu$ on $\C$ to be the Riesz measure for some
subharmonic function satisfying conditions (\ref{4}) and (\ref{5})
it is necessary and sufficient that the conditions (\ref{8}),
(\ref{13}), (\ref{22}), and
\begin{equation}\label{27}
    \sup\limits_{x\in\R} \int\limits_{1}^{\infty} \frac{\mu(B(0,t)) -
    \mu(B(x,t))}{t}dt < \infty.
\end{equation}
be fulfilled.
\end{theor}

\bproof

If a subharmonic function $v$ satisfies (\ref{4}) and (\ref{5}),
then, by Theorem 2, its Riesz measure $\mu$ satisfies (\ref{8})
and the measures $\mu_t$ converge to the measure $\mu_{\infty} =
(\s_+ +\s_-)(2\pi)^{-1}m_1$. The last measure does not charge any
circle $|w|=R$, therefore (\ref{7}) implies that $\mu (B(0,R)) = C
R + \overline{o}(R)$ as $R\to\infty$. Hence we get (\ref{13}) and
(\ref{22}). By theorem 1, $v(z)=V(z)+A_1y+A_2$. So the function
$V$ is also bounded from above on $\R$. Since Theorem FL, we see
that the same is true for the function $\tilde V$ from (\ref{28}).
Now Lemma 3 implies (\ref{27}).

Conversely, if a measure $\mu$ satisfies (\ref{8}), (\ref{13}),
and (\ref{22}), then, by Lemma 3, the subharmonic function $V$ has
the Riesz measure $\mu$ and satisfies (\ref{4}) and (\ref{5}).
Theorem is proved.

\begin{theor}.
For a measure $\mu$ on $\C$ to be the Riesz measure for some
subharmonic function $v$ with the property (\ref{4}) such that the
family of translations $\{v(z+h)\}_{h\in\R}$ is a relatively
compact subset in $\cal D'(\C)$, it is necessary and sufficient
that the conditions (\ref{8}), (\ref{13}), (\ref{22}), and
\begin{equation}\label{30}
    \sup_{x\in\R} \left|\int\limits_{1}^{\infty} \frac{\mu(B(0,t)) -
    \mu(B(x,t))}{t}dt\right|< \infty
\end{equation}
be fulfilled.
\end{theor}

\bproof. If a function $v$ satisfies (\ref{4}), and the function
$\tilde v$ from (\ref{31a}) is bounded  on $\R$, then $v$ is
bounded from above on $\R$. Theorem 3 implies conditions
(\ref{8}), (\ref{13}), and (\ref{22}) for the Riesz measure $\mu$
of the function $v$. By Theorem 1, $\tilde v$ equals $\tilde V$
from (\ref{28}) on $\R$ up to a constant term, therefore Lemmas 3
and 4 imply (\ref{30}).

Conversely, if a measure $\mu$ satisfies (\ref{8}), (\ref{13}),
(\ref{22}), and (\ref{30}), then $V$ from (\ref{23}) satisfies
(\ref{4}), and the function $\tilde V$ is bounded on $\R$.

Hence the assertion of Theorem 4 follows from Lemma 4.
\eproof

\begin{center}
{\bf \S 3}
\end{center}
\bigskip

A continuous function $F(z)$ on a closed strip $\{z=x+iy:\
x\in\R,\ |y|\le H\}$ with $H\ge 0$ is {\it almost periodic} if the
family of translations $\{F(z+h)\}_{h\in\R}$ is a relatively
compact set with respect to the topology of uniform convergence on
the strip; a function is almost periodic on an open strip (in
particular, on $\C$), if it is almost periodic on every closed
substrip of a finite width.

A measure (maybe complex) $\mu$ on $\C$ is called {\it almost
periodic} if for any test-function $\f \in \cal D(\C)$ the
convolution $\int \f(w+t) d\mu(w)$ is an almost periodic function
in $t\in\R$ (\cite{r}).

The following statement is valid:

{\bf Theorem R} (Theorem 1.8  \cite{r}). {\it For a measure $\mu$
to be almost periodic it is necessary and sufficient that the
following condition be fulfilled: for each sequence
$\{h_n\}\subset\R$ there exists a subsequence $\{h_n'\}$ such that
the convolutions $\int \f(w+x+h_n') d\mu(w)$ converge uniformly
with respect to  $x\in\R$ and functions $\f\in L$, where $L$ being
a compact subset in $\cal D(\C)$. Moreover, for a measure to be
almost periodic it is sufficient to check this condition only for
all single-point sets $L\subset \cal D(\C)$}.

If we take $L=\{\f(z+iy)\}_{|y|\le H}$ for some $\f\in \cal
D(\C)$, we obtain that  the convolutions $\int \f(w+z+h_n')
d\mu(w)$ actually converge uniformly on any strip $|y|\le H$,
hence the function $\int \f(w+z) d\mu(w)$ is almost periodic on
$\C$.

Further, a subharmonic function $v$ on $\C$ is called {\it almost
periodic}, if the measure $v(z)dm_2(z)$ is almost periodic (see
\cite{rrf};  equivalent definition see \cite{fr}).

It follows from definition that the Riesz measure of an almost
periodic subharmonic function is also  almost periodic.
Conversely, each almost periodic measure is the Riesz measure of
some almost periodic subharmonic function (see \cite{rrf}, where
to be investigated the case of a strip as well). Note  that the
family of translations $\{v(z+h)\}_{h\in\R}$ is a relatively
compact subset of $\cal D'(\C)$ for every almost periodic
subharmonic function $v$ on $\C$. Note also that each almost
periodic subharmonic function is bounded from above on every
horizontal strip of a finite width (see \cite{rrf}).

\bigskip

Here we obtain the following result:

\begin{theor}.
A necessary and sufficient conditions for a measure $\mu$ on $\C$
to be the Riesz measure of some almost periodic subharmonic
function at most of linear growth is that the measure be almost
periodic and satisfied (\ref{13}), (\ref{22}), (\ref{8}), and
(\ref{30}).
\end{theor}

The proof of the theorem bases on the following lemmas.

\begin{lem}. Suppose subharmonic functions $v_1(z)$ and
$v_2(z)$ on $\C$ with the common Riesz measure satisfy (\ref{4})
and (\ref{5}); then $v_1(z) = v_2(z) + p_1 + p_2 y$. Further,
 if
$$
\sup\limits_{\R}v_1(x) =\sup\limits_{\R}v_2(x) \quad \hbox{and}
\quad \s_+(v_1) = \s_+(v_2),
$$
where $\s_+$ be defined in (\ref{11}), then $v_1\equiv v_2$.
\end{lem}

\bproof The first part follows from Theorem 1, the second one is
evident.  \eproof

\begin{lem}. Let a subharmonic on $\C$ function $v$ be
satisfied (\ref{4}) and a family $\{v(z+h)\}_{h\in\R}$ be a
relatively compact subset of $\cal D'(\C)$;  if $v(z+h_n)\to
v^*(z)$ in the space $\cal D'(\C)$, then the family
$\{v^*(z+h)\}_{h\in\R}$ is a relatively compact subset as well and
\begin{equation}\label{31}
    \sup\limits_{x\in\R} v^*(x) \le \sup\limits_{x\in\R} v(x),
\end{equation}
\begin{equation}\label{32}
    \inf\limits_{t\in\R} \int\limits_{|z-t|<1} v^*(z)dm_2(z)
    \ge \inf\limits_{t\in\R} \int\limits_{|z-t|<1} v(z)dm_2(z),
\end{equation}
\begin{equation}\label{33}
    \s_+(v^*) \le \s_+(v),  \qquad     \s_-(v^*) \le \s_-(v).
\end{equation}
\end{lem}

\bproof\, Put $M=\sup\limits_{x\in\R}v(x)$. By Theorem FL, we get
for any $\e>0$
$$
    v(z)\le M +2\e \max \{\s_+(v),\s_-(v)\}, \quad |y|<2\e.
$$
Let $\f$ be a function from $\cal D(\C)$  such that $\f$ depends
only on $|z|$, $\f\ge 0$, $\f(z)=0$ for $|z|\ge\e$,
$\int\f(z)dm_2(z)=1$. Then
 $$
 (v\ast\f)(z)\le M +2\e \max \{\s_+(v),\s_-(v)\}, \quad |y|<\e.
 $$
Therefore,
 $$
 (v^*\ast\f)(z)\le M +2\e \max \{\s_+(v),\s_-(v)\}, \quad |y|<\e.
 $$
Note that $v^*$ is a subharmonic function, hence
$(v^*\ast\f)(z)\ge v^*(z)$. Since $\e$ is arbitrary, we obtain
(\ref{31}). By the same argument, for all $y\in\R$
 $$
 \sup_{x\in\R}v^*(x+iy)\le\sup_{x\in\R}v(x+iy).
 $$
Therefore we obtain (\ref{4}) and (\ref{33}).

Further, the functions $v(z+h_n)$ are integrable on every disc and
uniformly bounded from above, therefore we can replace the
convergence of measures $v(z+h_n)dm_2(z)$ in the sense of
distributions by the weak convergence of measures. Since the limit
measure $v^*(z)dm_2(z)$ does not charge any circle, we have
$$
\lim\limits_{n\to\infty} \int\limits_{|z-t|<1} v(z+h_n)dm_2 (z) =
\int\limits_{|z-t|<1} v^*(z)dm_2 (z)
$$
for each $t\in\C$. Hence we get (\ref{32}). Taking into account
Lemma 4, we see that the family $\{v^*(z+h)\}_{h\in\R}$ is a
compact subset of $\cal D'(\C)$. Lemma is proved.

\begin{lem}. Under the conditions of the previous lemma, suppose that the Riesz measure
$\mu$ of the function $v(z)$ is almost periodic. Then inequalities
(\ref{31}) - (\ref{33}) turn into equalities, the Riesz measure
$\mu^*$ of the function $v^*(z)$ becomes almost periodic, and
there is a subsequence $\{h_{n'}\}$ such that for every $\f \in
\cal D(\C)$
\begin{equation}\label{36}
    \lim\limits_{n'\to\infty} \sup\limits_{t\in\R}
    \left|\int \f(w-t-h_{n'})d\mu(w) - \int \f(w-t)d\mu^*(w) \right| =
    0.
\end{equation}
\end{lem}

\bproof

For all $\f \in \cal D(\C)$ we have
$$
    \lim\limits_{n\to\infty} \int \f(z-h_n)v(z)dm_2(z) =
    \lim\limits_{n\to\infty} \int \f(z)v(z+h_n)dm_2(z) = \int
    \f(z)v^*(z)dm_2(z).
$$
Since $\mu = (2\pi)^{-1} \triangle v$, we obtain
\begin{equation}\label{37}
    \lim\limits_{n\to\infty} \int \f(z-h_n)d\mu(z) = \int
    \f(z)d\mu^*(z).
\end{equation}

From Theorem R it follows that there is a subsequence $\{h_{n'}\}$
such that for any $\f\in \cal D(\C)$ the almost periodic functions
$\int \f (z-t-h_{n'})d\mu(z)$ converge to an almost periodic
function uniformly in $t\in\R$. If we replace $z$ by $z-t$ in
(\ref{37}), then we get (\ref{36}). Consequently, the function
$\int\f(w-t)d\mu^*(w)$ is almost periodic in $t\in\R$, and $\mu^*$
is an almost periodic measure.

Pass to a subsequence again if necessary, we may assume that the
functions $v^*(z-h_n)$ converge in the space $\cal D'(\C)$ to some
subharmonic function $v^{**}(z)$ with the Riesz measure
$\mu^{**}$. Therefore,
$$
    \lim\limits_{{n'}\to\infty} \int \f(z+h_{n'})d\mu^*(z) = \int
    \f(z)d\mu^{**}(z).
$$
On the other hand, it follows from (\ref{36}) that
$$
    \lim\limits_{n\to\infty}\left|\int \f(w)d\mu(w) - \int \f(w+h_n)d\mu^*(w) \right| = 0.
$$
Here $\f$ is an arbitrary function from $\cal D(\C)$, hence,
$\mu^{**} = \mu$. By Lemma 5, we get $v^{**}(z) = v(z) + D_1 + D_2
y$.

Since (\ref{31}) is valid for pairs $v, \, v^{*}$ and $v^{*},\,
v^{**}$, we get $D_1 \le 0$. Then (\ref{32}) for pairs $v, \,
v^{*}$ and $v^{*},\, v^{**}$ implies $D_1 \ge 0$, and we obtain
$D_1=0$ and the equality in (\ref{31}) and (\ref{32}). By the same
way, we obtain the equalities in (\ref{33}). Lemma is proved.

\bigskip

\noindent{\bf Proof of Theorem 5}. The necessity follows
immediately from Theorem 4. Let us prove a sufficiency. Suppose
$\mu$ satisfies the conditions of Theorem 5. Let $V$ be the
function from (\ref{23}), and let $\{h_n\}\subset\R$ be an
arbitrary sequence. It follows from Theorem 4 that the family
$\{V(z+h_n)\}$ is a relatively compact subset of $\cal D'(\C)$.
Therefore we can assume without loss of generality that
$V(z+h_n)\to V^*(z)$ in $\cal D'(\C)$. To prove the Theorem, we
need to check that
$$
\int \f(z-t-h_n)V(z) dm_2(z) \to \int \f(z-t)V^*(z) dm_2(z).
$$
uniformly in $t\in \R$ for any $\f\in \cal D(\C)$.

Assume the contrary. Then there is $\f_0\in \cal D(\C)$, $\e_0>0$,
and $\{t_n\}\subset\R$  such that
\begin{equation}\label{38}
    \left|\int \f_0(w)V(w+h_n+t_n)dm_2(w) - \int
    \f_0(w)V^*(w+t_n)dm_2(w) \right| \ge \e_0
\end{equation}
(if necessary we can replace the sequence $\{h_n\}$ by a
subsequence).

We may assume also that $V(z+h_n+t_n)\to V^{**}(z)$,
$V^*(z+t_n)\to V^{***}(z)$ in $\cal D'(\C)$ as $n\to\infty$. By
$\mu^*$, $\mu^{**}$, $\mu^{***}$ denote the Riesz measures of the
functions $V^*$, $V^{**}$, $V^{***}$, respectively. Then we have
\begin{equation}\label{39}
\lim\limits_{n\to\infty} \int \f(w-h_n-t_n)d\mu(w) = \int
\f(w)d\mu^{**}(w),
\end{equation}
\begin{equation}\label{40}
\lim\limits_{n\to\infty} \int \f(w-t_n)d\mu^*(w) = \int
\f(w)d\mu^{***}(w)
\end{equation}
for any $\f\in \cal D(\C)$.

On the other hand, the measure $\mu$ satisfies (\ref{36}). Hence
the integrals in left-hand sides of (\ref{39}) and (\ref{40}) have
the same limit, and $\mu^{**} = \mu^{***}$.

By Lemma 7 we obtain
$$
\sup\limits_{\R} V^{**}(x) = \sup\limits_{\R} V(x)
=\sup\limits_{\R} V^{*}(x) = \sup\limits_{\R} V^{***}(x)
$$
and
$$
\s_+(V^{**}) = \s_+(V) = \s_+(V^{*}) = \s_+(V^{***})
$$
Using Lemma 5, we get $V^{**} \equiv V^{***}$. This contradicts
(\ref{38}). Theorem 5 is proved.

\bigskip

Now let $d$ be a divisor in $\C$, i.e., a sequence
$\{a_k\}\subset\C$ without finite limit points such that each
value may appear with a finite multiplicity. A divisor is called
{\it almost periodic} if the discrete measure supported at the
points $a_k$ with mass at each point equals the multiplicity of
the point in the sequence is almost periodic (see \cite{r1},
\cite{rrf}; in \cite{rrf} there is an equivalent geometric
definition). Moreover, almost periodic divisors are just the
divisors of entire functions\footnote{A divisor $\{a_k\}$ is the
divisor of a holomorphic function $f$ when zeros of $f$ coincide
with the values $\{a_k\}$ and a multiplicity of every zero equals
the multiplicity of corresponding $a_k$.} with almost periodic
modulus in every substrip $\{z=x+iy:\,|y|<H\}$. (\cite{rrf}).

\begin{theor}.
For a divisor $\{a_k\}$ to be the divisor of an entire function of
exponential type with the almost periodic modulus, it is necessary
and sufficient the following conditions be fulfilled:

a) The divisor is almost periodic,

b) there is a finite limit
$$
\lim\limits_{R\to\infty} \sum\limits_{|a_k|<R} \frac{1}{a_k},
$$

c) $n(0,t) = O(t)$,

d) $n(0,t+1) - n(0,t) = \overline{o}(t)$,

e)$$
    \sup\limits_{x\in\R} \left| \int\limits_{1}^{\infty} \frac{n(0,t) -
    n(x,t)}{t}dt \right| < \infty;
$$
here $n(c,t) = \card \{ k: \; |a_k - c|\le t \}$.

\end{theor}

\bproof By Theorem 6, conditions a) -- e) mean just the existence
of an almost periodic subharmonic function $v$ at most of linear
growth with the Riesz measure supported at the points $\{a_k\}$
with mass equals a multiplicity of the point in the sequence. Then
$v(z) = \log |f(z)|$ for an entire function $f$ of exponential
type such that the divisor of $f$ is $\{a_k\}$. Since $v$ is an
almost periodic, we get that $|f|$ is almost periodic (see
\cite{rrf}). Theorem is proved.

\bigskip

Now we consider the periodic case.

\begin{theor}.
A necessary and sufficient condition for a measure $\mu$ on $\C$
to be the Riesz measure of some periodic subharmonic function with
period $1$ at most of linear growth is the measure be stable with
respect to the translation on $1$ and
\begin{equation}\label{44}
    \mu \{ z=x+iy: \; 0 \le x < 1, \; y\in \R \} < \infty.
\end{equation}

\end{theor}

\bproof Let $v$ be a subharmonic function such that $v(z+1) =
v(z)$. It is clear that its Riesz measure is stable with respect
to the translation on $1$. By Theorem 1, it follows that $\mu$
satisfies (\ref{13}). Using the equality
\begin{equation}\label{45}
    \mu \{ z=x+iy: \; 0\le x < 1, \; |y|<n\} = \frac{1}{2n}
       \mu \{ z=x+iy: \; -n\le x < n, \; |y|<n\},
\end{equation}
we get (\ref{44}).

Conversely, let $\mu$ be stable with respect to the translation on
$1$ and be satisfied (\ref{44}). Using (\ref{45}), we obtain
(\ref{13}). Then for any $r>0$

\begin{equation}\label{46}
 \mu \{ z: \; r \le |z| < r+1 \} \le \sum\limits_{n\in\Z,\; |n|\le r+1}
 \mu \{ z: \; x\in [0,1), \; r \le |z+n| < r+1 \}
\end{equation}

Fix $\d \in (0,1)$. For $|n| < r(1-\d)-1$ we have
$$
\{ z: \; x\in [0,1), \; r \le |z+n| < r+1 \} \subset \{ z: \; x\in
[0,1), \; |y| > r\d \}
$$
Hence (\ref{46}) is majorized by
\begin{equation}\label{46a}
2r(1-\d) \mu \{ z: \; x\in [0,1), \; |y| > r\d \} + (2\d r +5) \mu
\{ z: \; x\in [0,1), \; y \in \R \}.
\end{equation}
It follows from (\ref{44}) that for any $\e > 0$  there exist $\d
> 0$ and $r_0<\infty$ such that for $r\ge r_0$ (\ref{46a}) is less than $r\e$.
This yields (\ref{22}).

Further, take $R>r>1$. We have
$$
\left|\; \int\limits_{r<|z|<R} \frac{d\mu (z)}{z} -
\int\limits_{r<|[x]+iy|<R} \frac{d\mu (z)}{z}
    \right|\le
$$
\begin{equation}\label{47}
    \le \frac{1}{r-1} \mu \{ z:\; r-1 <|z|<r+1\}
     +  \frac{1}{R-1} \mu \{ z:\; R-1 <|z|<R+1\};
\end{equation}
where $[x]$ being the integral part of a real $x$. By (\ref{22}),
the right-hand side of (\ref{47}) tends to zero as $r\to\infty$.
Then we obtain
$$
    \int\limits_{r<|[x]+iy|<R} \frac{d\mu (z)}{z} =
    \sum\limits_{n\in\Z:\;|n|\le R}\int\limits_{x\in[0,1),\; r<|n+iy|<R} \frac{d\mu
    (z)}{z+n}=
$$
\begin{equation}\label{48}
    = \int\limits_{x\in[0,1),\; y\in \R} \sum\limits_{n\in\Z:\;r<|n+iy|<R}
    \frac{\overline{z} + n}{|z+n|^2}d\mu(z)
\end{equation}
Now for any $x\in [0,1)$, $y\in\R$ we have
$$
\sum\limits_{n\in\Z:\;|n+iy|>r} \frac{|\overline{z}|}{|z+n|^2} \le
\sum\limits_{n\in\N\bigcup \{0\}:\;n^2 > r^2 -y^2}
\frac{1+|y|}{n^2 + y^2} + \sum\limits_{n\in\N:\;n^2 > r^2 -y^2}
\frac{1+|y|}{(n-1)^2 + y^2} \le
$$
\begin{equation}\label{49}
\le 2 \int\limits_{\sqrt{\max \{0,r^2-y^2\}}}^{\infty}
\frac{1+|y|}{t^2 + y^2}dt + \frac{2(1+|y|)}{y^2} \chi_r(y);
\end{equation}
here $\chi_r(y)$ is a characteristic function of interval
$(\sqrt{r^2-1},\infty)$. Besides,
$$
- \sum\limits_{n\in\Z:\; r<|n+iy|<R} \frac{n}{|z+n|^2} =
\sum\limits_{n\in\N:\; r < |n+iy|<R} n \left( \frac{1}{|z-n|^2} -
\frac{1}{|z+n|^2} \right) =
$$
\begin{equation}\label{50}
= \sum\limits_{n\in\N:\; r^2 < n^2+y^2<R^2} \frac{2n^2
x}{((n-x)^2+y^2)((n+x)^2+y^2)}
\end{equation}

It is easy to see that the right-hand side of (\ref{50}) is also
majorized by (\ref{49}). Both terms monotonically decrease to 0 as
$r\to\infty$, hence (\ref{48}) tends to 0 as $r\to\infty$
uniformly in $R$. It follows from (\ref{47}) and (\ref{48}) that
(\ref{8}) is valid. Finally, by (\ref{28}) and (\ref{24}), the
integral
$$
\int\limits_{1}^{\infty} \frac{\mu(B(0,t)) - \mu (R(z,t))}{t} dt
$$
is bounded for $z=x\in [0,1]$. Since $\mu$ is stable with respect
to the translation on $1$, we get (\ref{30}). Now the assertion of
Theorem 7 follows from Theorem 5.

\eproof

\bigskip

{\bf Consequence.} {\it A necessary and sufficient conditions for
a divisor $\{a_k\}$ to be the divisor of an entire periodic (with
period 1) function of exponential type bounded on real axis is
that the divisor be periodic with period 1 and its restriction to
the strip $\{z:\; 0\le x < 1,\,y\in\R\}$ be finite.}

In this case the corresponding  function is a finite product of
elementary functions $\sqrt{1-\cos 2 \pi (z-\gamma_k)}$ with $\Re
\gamma_k \in [0,1)$; it is unique up to a multiplier $Ce^{i\nu
z}$, $\nu\in\R$.

\bigskip

Department of Mathematics, Kharkov National University,

Svobody sq.,4, Kharkov 61077, Ukraine,

Sergey.Ju.Favorov@univer.kharkov.ua

\end{document}